# Design and analysis of demand-adapted railway timetables

David Canca[a]*, Eva Barrena[a], Encarnación Algaba[a], Alejandro Zarzo[b]

[a]*E.T.S. Ingenieros,Universidad de Sevilla, Camino de los Descubrimientos s/n, 41092. Seville, Spain*
[b]*E.T.S.I.I. Universidad Politécnica de Madrid, José Gutiérrez Abascal, 2, 28006. Madrid, Spain*

**Abstract**

Railway scheduling and timetabling are common stages in the classical hierarchical railway planning process and they perhaps represent the step with major influence on user's perception about quality of service. This aspect, in conjunction with their contribution to service profitability, makes them a widely studied topic in the literature, where nowadays many efforts are focused on improving the solving methods of the corresponding optimization problems. However, literature about models considering detailed descriptions of passenger demand is sparse. This paper tackles the problem of timetable determination by means of building and solving a non-linear integer programming model which fits the arrival and departure train times to a dynamic behavior of demand. The optimization model results are then used for computing several measures to characterize the quality of the obtained timetables considering jointly both user and company points of view. Some aspects are discussed, including the influence of train capacity and the validity of Random Incidence Theorem. An application to the C5 line of Madrid rapid transit system is presented. Different measures are analyzed in order to improve the insight into the proposed model and analyze in advance the influence of different objectives on the resulting timetable.

## 1. Introduction

Railway scheduling and timetabling are studied at the third stage in a hierarchical schema composed by five consecutive steps: Analysis of demand, line planning, scheduling and timetabling, rolling-stock and crew management (Bussieck et al., 1997, Bussieck, 1998). Timetabling may be one of the topics with major influence on users' quality of service perception. In fact, this subject has received the attention of many researches and has been widely studied in the literature, where, nowadays the main lines of research tend to improve the solving methods of the corresponding optimization problems. Cacchiani, 2008, presented a review of optimization models arising railway planning and centered on the role of scheduling and timetabling inside the classic hierarchical planning approach. Timetabling has specific characteristics depending on the type of service to be provided. So, if the frequency of the service is high enough, to specify this frequency becomes more important than the detailed hourly schedule itself, e.g., for rapid transit systems information about headway between trains is more important than knowing the exact departure/arrival times (see Higgins and Kozan, 1998).

In the train timetabling phase, all arrival and departure times are obtained. This task is frequently done subject to the periodicity of the system. Periodic scheduling, initially proposed by Voorhoeve, 1993, based on the formulation of Serafini and Ukovich, 1989, has been followed by other authors like Natchtigall, 1994, Natchtigall and Voget, 1997, Liebchen and Möhring, 2002, Liebchen and Peeters, 2002, 2009, and Chierici et al., 2004. Periodic timetables have some advantages in the case of Metropolitan railways since a simpler service is allowed. In fact, regular timetables are easily memorized by users and also can be computed with less effort (Wardman et al., 2004). On the contrary, in the case of middle and long distance networks, the demand is usually captive to the schedules (see Pham 2004, Vansteenwegen and van Oudheusden, 2006). In this kind of scenarios, when passenger demand seems to be dependent on railway scheduling, the periodic regular approach may not be efficient. Even in the situation of subway and local railways, when a detailed description of passenger demand is considered as one of the main factors to determine timetables, periodic timetables are suboptimal from a user's perspective.

Note that an inefficient design of timetables could lead to mode change decisions of users, choosing alternative transportation modes (with the consequent drop of network profitability). Non-periodic timetabling is especially important in long corridors with high traffic densities. This approach allows for obtaining optimal departure and arrival times after the line planning step (see Kwan and Mistry, 2003, Caprara et al., 2006 or Ingolotti et al., 2006). It should be mentioned that most of the existent approaches do not make realistic hypothesis about the behavior of the demand of each line over a full day of operation. In fact, as we mention above, the analysis of the demand is commonly done only in the first stage of the planning process, generating a peak hour origin-destination matrix that is used to determine frequencies in the line planning phase. Interesting exceptions are the works of Chierici et al., 2004, and the posterior extension published by Cordone and Redaelli, 2011, when they relax the common assumption of regular demand to consider the influence between timetable quality and captured demand. In this context, once the timetable is obtained, some trains are removed to accomplish the transport demand for each non-peak hour. In any case, few authors deal with the problem of relating the scheduling of the units to the quality of service (Natchtigall and Voget, 1997, Vansteenwegen and van Oudheusden, 2006) and the actual capacity of the network (see Burdett and Kozan, 2006 for an analytical approach, Abril et al., 2008, Canca, 2009, Canca et al., 2011 for several optimization approaches and Barber et al., 2007, for a description of simulation based systems). We would like to emphasize that our approach to timetables design is subsequent to the line planning phase, so that the lines are already defined and no more modified by the daily variations in the demand. Moreover, the timetable model proposed here adapts the solution for each line to realistic descriptions of these variations.

The main goal of this paper is to analyze the performance of railway timetables by means of some proposed relative measures of quality of service and line profitability which are explicitly constructed and afterward studied in detail. For that aim, an improved and more accurate version of the approach already described in Canca et al., 2011 is considered in this work by applying two consecutive linked models (see Figure 1).

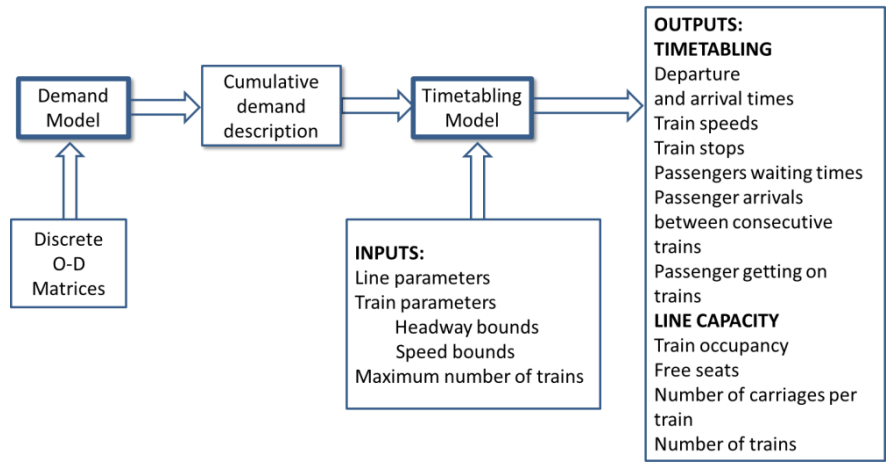

Figure 1. Scheme of the modelling framework

The first one provides a continuous description of the mobility demand (instead of using the more conventional discrete peak-hour origin-destination matrix picture) by fitting cumulative demand data with certain sums of sigmoidal functions (see Canca et al., 2011 and Section 2 below), being able to simulate the demand peaks along the planning horizon in a very accurate way. Here, with respect to the one used in Canca et al., 2011, this demand model has been updated to ensure a good approximation outside the planning horizon, an important aspect to be considered in the second model for the required actual computation of arrivals when certain trains are not launched.

Some characteristics of the railway network together with the cumulative demand description provided by this first model, will serve as the required inputs for the second model, which is the cornerstone of our approach. It is a non-linear integer programming model designed to obtain optimal timetables using different aspects of the network and demand behavior (see Canca et al., 2011 and Section 3 below). In comparison with the preliminary version, used in Canca et al., 2011, here this second model has been redefined to extend and improve many aspects. Changes and modifications included in this new version have been mainly aimed at the improvement of three fundamental aspects of the model: its realism, its efficiency when evaluating quality of service measures of the network and its accurateness when analyzing profitability measures of the network.

Thus, the inclusion of train launching decisions as well as variable stop times makes the model much more flexible and realistic at the cost of increasing the computational time since it becomes necessary to introduce a complete new set of binary variables to model the launching process. To diminish this undesired consequence, a new ordering constraints set is introduced to reduce the search space in ordering decisions.

Stop times are now partially dependent on the demand, so that they are obtained considering reasonable bounds, passenger arrivals between consecutive trains and free capacity at each stop. All of these new characteristics, together with the possibility of taking non stopping decisions, give rise to the possibility of a much more efficient evaluation of quality of service measures. Moreover, when analyzing profitability, a key aspect is to select the most convenient capacity. For this selection (see Section 3.3, below), the use of load factors is now proposed to compute performance measures from a service provider point of view. Also, new capacity curves considering constant number of carriages are included to discuss the appropriate capacity selection when global capacity can serve all the demand in conditions of low occupancy.

Finally, our new models and methodology have been tested (see Section 4) in a real scenario which is a piece of the C5 line belonging to the Madrid suburban railway. A final analysis of the profitability of the offered services for each line requires a global study of the network. Some aspects like maintenance, rolling stock and crew management should be taken into account (Lindner and Zimmermann, 2005). The complexity of the global problem leads to the search for certain measures which would help to select an appropriate timetable, subject to a global economic post-study, in charge of the service provider.

## 2. Passenger demand analysis

Usually, the demand mobility is characterized by means of the so-called *Origin-Destination* matrices. These kinds of matrices are commonly used in transportation planning analysis, where planners work with different matrices for a design day, each one for a defined time interval, e.g., the peak hour matrix or the total daily demand matrix. The elements of these matrices are usually obtained by means of more or less rigorous extrapolation techniques applied to the data given by mobility surveys (of course in a discrete way). In this paper, we follow a different approach considering a generic continuous representation of the mobility demand $OD(t)$ along the daily planning horizon. So, $OD(t)$ is a square matrix as many rows and columns as stations in the network. Each entry, denoted by $f_{ij}(t)$ represents the continuous daily evolution of the travel demand from station $i$ to station $j$.

Although these demand curves can quite differ from a problem to another, a common and relevant characteristic of all of them is the existence of certain demand peaks (local maxima) in certain instants of time along the day. These peaks are associated to rush hours and generally reduced to two or three per day. Following the ideas already described in Canca et al., 2011, in our approach to the passengers demand, we consider the cumulative or aggregate demand $F_{ij}(t)$, given by,

$$F_{ij}(t) = \int_o^t f_{ij}(s)\,ds, \qquad (0)$$

thus, the number of passengers arriving at the $i$-th station bound to the $j$-th during the time interval $[t_1, t_2]$ is given by,

$$F_{ij}^{[t_1,t_2]} = F_{ij}(t_2) - F_{ij}(t_1) = \int_o^{t_2} f_{ij}(s)\,ds - \int_o^{t_1} f_{ij}(s)\,ds = \int_{t_1}^{t_2} f_{ij}(s)\,ds, \qquad (0)$$

as it is shown in Figure 2. In this figure, it can be observed how the demand changes over time.

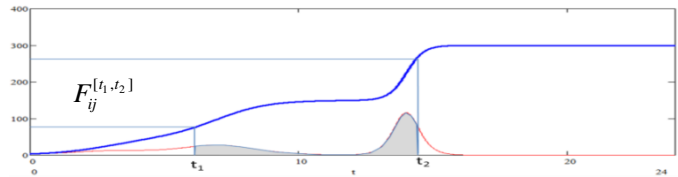

Figure 2. Demand and cumulative demand functions

As the demand will be used to obtain train departure times and provide measures to characterize the performance of timetables, an analytical expression of it would be of interest. Taking into account the shape of the cumulative demand between each pair of stations, it makes sense to consider an approximation given by a linear combination of a variable number, $M$, of sigmoid curves, i.e.,

$$F_{ij}(t) = \sum_{r=1}^{M} \frac{K_{ij}^{r}}{1+e^{-\beta_{ij}^{r}(t-x_{ij}^{r})}}. \tag{0}$$

The number of terms used in the approximation, i.e., $M$, symbolizes the maximum number of demand peaks along the day. In each sigmoid function, the parameters $K_{ij}^{r}$, $x_{ij}^{r}$ and $\beta_{ij}^{r}$ stand for the asymptotic value, the deviation with respect to the time and the slope, respectively.

This approximation is fully characterized by a number of parameters (including $M$) which are determined by solving a set of appropriate least squares minimization problems given by,

$$Min\left(\sum_{t=1}^{T}(y_{ij}(t)-F_{ij}(t))^2\right),$$

$s.t.$

$$F_{ij}(t) = \sum_{r=1}^{M} \frac{K_{ij}^{r}}{1+e^{-\beta_{ij}^{r}(t-x_{ij}^{r})}},$$

$$\sum_{r=1}^{M} K_{ij}^{r} \leq \max_{t\leq T} y_{ij}(t) = y_{ij}(T),$$

$$K_{ij}^{r},\ \beta_{ij}^{r},\ x_{ij}^{r} \geq 0. \tag{0}$$

In other words, the dynamic adjust of these functions is made by minimizing the sum of quadratic errors as the difference between the proposed demand function $F_{ij}(t)$ and the demand data $y_{ij}(t)$ attained, in a standard day, for each pair of stations $i, j, i \neq j$, on certain instants depending on data availability.

It should be mentioned that real demand data can be given in form of a set of origin-destination matrices corresponding to different intervals along the planning horizon, i.e., configuring a discrete description of mobility demand. Also, entrance and exit counts at stations (in case of availability) can be used to build the cumulative demand function between each pair of stations. In our knowledge, nor timetabling neither scheduling models have been previously proposed considering variable demand data as input to adjust departure and arrival times in a long period.

### 3. Timetabling optimization model

In this section, we present an enhanced optimization model based on a previous work (Canca et al., 2011) in order to determine the train departure times from the origin station as well as the arrival/departure times at/from each one of the stations, for each train, along a line. As mentioned before, this approach, unlike other models in the existent literature, is based on the cumulative demand approximations described in the above section and, as consequence, a set of specific variables and constraints are used to model passenger behavior. We emphasize that this model can be applied to other kind of characterizations of the demand, even it can be adaptable for its application to problems in which discrete demand functions are used, as the ones typically obtained from mobility surveys. The aim of this formulation is, not only the calculation of train arrival/departure times, but also obtaining measures that allow for the analysis of timetables quality.

In order to describe the optimization model, in the following subsections, the notation, constraints and a discussion on objective function and quality measures will be introduced.

### 3.1. Notation

Given a generic line, the notation of sets, parameters and variables, used in the model, is given in Table 1.

| *Sets and parameters* | |
|---|---|
| Symbol | Definition |
| $H := \{1, 2, ..., S\}$ | Set of stations. |
| $N := \{1, 2, ..., K\}$ | Set of trains. See equation (9) for an explanation of the value of $K$. |
| $L := \{1, 2, ..., S-1\}$ | Set of segments between each pair of stations. |
| $CAP_{min}$ | Minimum capacity of every train unit., i.e. sum of capacities of the minimum carriage number allowed. |
| $T$ | Planning horizon, usually one day (time units). |
| $lng(i,i+1)$ | Length of the segment corresponding to stations $i$ and $i+1$. |
| $Mint_{stop}$, $Maxt_{stop}$ | Lower and upper bounds for stop times at stations. |
| $t_{saf}(i,k)$ | Safe headway time after departure of train $k$ from station i. |
| $t_{acc}(i,i+1)$ | Acceleration time needed to reach pure running speed after station $i$. |
| $t_{dec}(i,i+1)$ | Deceleration time needed to stop from running speed before station $i+1$. |
| $g$ | Passenger flow rate, pax/min. Door open and close times are considered negligible. |
| *Variables* | |
| Symbol | Definition |
| $t_i^k$ | Departure time of train $k$ at the $i$-th station. |
| $\overline{V}_k(i,i+1)$ | Unit travel time of train $k$ in the segment between stations $i$ and $i+1$. |
| $FS_i^k$ | Available capacity in the $k$-th train when it leaves the $i$-th station. |
| $CAP_k$ | Capacity of train $k$. |
| $t_{stop}(i,k)$ | Stop time of train $k$ at station $i$. |
| $\delta^k$ | Binary variable which takes value 1 if train with order $k$ is programmed. |
| $\beta_i^k$ | Binary variable which takes value 1 if train $k$ stops at station $i$. |
| $y^k$ | Interdeparture time between train $k$ and $k+1$ at the origin station. |
| *Passenger arrival variables* | |
| Symbol | Definition |
| $N_{ij}^k = F_{ij}^{[t_i^{k-1}, t_i^k]}$ | Number of people who arrive at the $i$-th station with destination to the $j$-th station during the time interval $[t_i^{k-1}, t_i^k]$. |
| $NAD_i^k$ | Number of people who arrive at the $i$-th station during the time interval $[t_i^{k-1}, t_i^k]$. |
| $NAO_j^k$ | Number of people who arrive at any station $i$ ($i \leq j$) with destination to the $j$-th station during the time interval $[t_i^{k-1}, t_i^k]$. |
| $ns_{ij}^k$ | Number of people who arrive at the $i$-th station with destination to the $j$-th station during the time interval $[t_i^{k-1}, t_i^k]$ and get on the $k$-th train. |
| $ne_{ij}^k$ | Number of people who arrive at the $i$-th station with destination to the $j$-th station during the time interval $[t_i^{k-1}, t_i^k]$ and do not get on the $k$-th train. |
| *Passenger waiting variables* | |
| Symbol | Definition |
| $E_{ij}^k$ | Number of people who arrive at the *i*-th station with destination to the *j*-th station before $t_i^{k-1}$ (i.e., during the interval $[0, t_i^{k-1}]$) and are waiting on platform before the departure of the *k*-th train. |
| $EAD_i^k$ | Number of people who arrive at the *i*-th station before $t_i^{k-1}$ (i.e., during the interval $[0, t_i^{k-1}]$) and are waiting on platform before the departure of the *k*-th train. |
| $EAO_j^k$ | Number of people who arrive at any station *i* ($i \leq j$) before $t_i^{k-1}$ (i.e., during the interval $[0, t_i^{k-1}]$) and are waiting on platform before the departure of the *k*-th train with destination to the *j*-th station. |
| $ee_{ij}^k$ | Number of people who arrive at the *i*-th station with destination to the *j*-th one, before $t_i^{k-1}$ (i.e., during the interval $[0, t_i^{k-1}]$) and do not get on the *k*-th train. |
| $es_{ij}^k$ | Number of people who arrive at the *i*-th station with destination to the *j*-th one, before $t_i^{k-1}$ (i.e., during the interval $[0, t_i^{k-1}]$) and get on the *k*-th train. |

***Passenger departure variables***

| Symbol | Definition |
|---|---|
| $S_{ij}^k$ | Number of people who get on the *k*-th train at the *i*-th station with destination to the *j*-th one. |
| $SAD_i^k$ | Number of people who get on the *k*-th train at the *i*-th station. |
| $SAO_j^k$ | Number of people who get on the *k*-th train with destination to the *j*-th station. |

Table 1. Symbols and definitions

We assume known the length and the minimum and maximum speed limitations in every segment, $l \in L$.

### 3.2. Constraints

We first present the constraints related to timetable variables. Equation (5) shows the relationship between departure times at consecutive stations along the line. An interesting extension, as proposed by Zhou and Zhong, 2005, considering acceleration and deceleration times, can be easily incorporated to this model by adding the term $(t_{acc} + t_{dec})\beta_{i+1}^k$ to the left side of constraints given in (5). In this case, $\overline{V}_k(i, i+1)$ refers to the "pure running" speed (see Zhou and Zhong, 2005). Notice that here we use the inverse unit travel time (inverse of average speed) in order to maintain linear expressions. Equation (6) defines the minimum headway between consecutive trains *k* and *k*+1, the arrival time of train *k*+1 is constrained according to the safety specifications for each train and station (see also Equation (7)).

$$t_{i+1}^k = t_i^k + t_{stop}(i+1, k) + lng(i, i+1)\overline{V}_k(i, i+1), \quad \forall k \in N, \ i \in H \setminus \{S\}. \qquad (0)$$

$$t_i^{k+1} - t_{stop}(i, k+1) \geq t_i^k + t_{saf}(i, k), \quad \forall k \in N \setminus \{K\}, \ i \in H. \qquad (0)$$

$$Mint_{stop} \ \beta_i^k \leq t_{stop}(i, k) \leq Maxt_{stop} \ \beta_i^k, \quad \forall k \in N, \ i \in H. \qquad (0)$$

$$\beta_1^k = 1, \ \beta_S^k = 1, \quad \forall k \in N. \qquad (0)$$

Also, train launching is controlled by binary variables $\delta^k$, if $\delta^k = 1$, train *k* is considered *active*, i.e., train *k* is launched in *[0,T)*, otherwise train *k* is inactive, see Figure 3. Obviously a pessimistic upper bound of train number is needed and this number is related to total demand and the minimum train capacity.

$$K = \frac{1}{CAP_{min}} \max_{i<S} \left\{ \sum_{p<i}^{r \geq i+1} F_{pr}^{[0,T]} \right\}. \quad (0)$$

Now, since train characteristics (capacity, speed, stopping) are incorporated in the model, without loss of generality, the train activation variables can be ordered. In this way, the model will select trains that are represented by lower values of the index $k$. These ordering constraints do not affect the optimal solution but contribute to reduce significantly the problem complexity.

$$\delta^k \geq \delta^{k+1}, \quad k \in N \setminus \{K\}. \quad (0)$$

When a train is *inactive*, the model forces departure time from station 1 beside $2T$, where cumulative demand function is approximately constant, i.e., arrivals between two consecutive trains are approximately zero. In general, any integer number $n \geq 2$ can be used.

$$2T(1-\delta^k) \leq t_1^k \leq T + nT(1-\delta^k), \quad \forall k \in N. \quad (0)$$

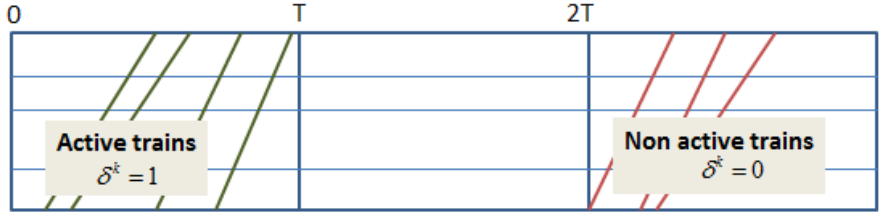

Figure 3. *Active* and *inactive* trains

Equations (5) and (6) consider departure times at stations and obviously these departure times should be directly related to passenger demand. Equation (12) implements this relationship, where passenger's arrival between the departure time of each pair of consecutive trains, $k$ and $k+1$, for every pair of stations, $i$ and $j$, is expressed using the approximation of cumulative demand described in Section 2.

$$N_{ij}^k = F_{ij}^{[t_i^{k-1}, t_i^k]} = \sum_{r=1}^{M} \frac{K_{ij}^r}{1+e^{-\beta_{ij}^r(t_i^k - x_{ij}^r)}} - \sum_{r=1}^{M} \frac{K_{ij}^r}{1+e^{-\beta_{ij}^r(t_i^{k-1} - x_{ij}^r)}}, \quad \forall i,j \in H,\ i<j,\ k \in N. \quad (0)$$

Note that arrivals between two consecutive trains, $k$ and $k+1$, are related to departure times $t_i^k$ and $t_i^{k+1}$ at the same station $i$. These set of equations are responsible for the presence of non-linearities in the model.

Equations (13) and (14) define aggregate arrival variables. The first one represents the sum of arrivals for passengers coming to station $i$ with destination to station $j$ ($j>i$). The second one computes the sum of passenger arrivals from several stations ($i<j$) to a certain one, denoted by index $j$. Without losing generality, we assume that the equations will be written enumerating the stations in a consecutive order. In the reverse way, it is sufficient to express the sum for $j < i$, if the same ordinals are used.

$$NAD_i^k = \sum_{\{j \in H\ :\ j>i\}} N_{ij}^k, \quad \forall i \in H \setminus \{S\},\ k \in N. \quad (0)$$

$$NAO_j^k = \sum_{\{i \in H\ :\ i<j\}} N_{ij}^k, \quad \forall j \in H \setminus \{1\},\ k \in N. \quad (0)$$

On the other hand, with the objective of calculating non-served passenger demand, the number of people who have arrived is decomposed into the sum of those passengers who will manage to get on the k-th train $\left(ns_{ij}^k\right)$ and those who will not achieve to do it $\left(ne_{ij}^k\right)$. Of course, this situation can occur when train capacity is lower than required to meet demand. Capacity is an important factor to determine quality of service and also affects line profitability because it is directly related to carriage number and its corresponding costs. In this model, train capacity is used as a variable, in contrast with other ones where capacity is considered as an input obtained after the line planning stage. Obviously, here, capacity can be fixed as well, if desired, so this is a more general approach. Notice that all posterior stations to the *i*-th station are considered as possible destination stations. Therefore,

$$N_{ij}^k = ns_{ij}^k + ne_{ij}^k, \quad \forall i,j \in H, \ i < j, \ k \in N. \tag{0}$$

In the same way, the number of passengers who wait for train *k*, having arrived at the platform before $t_i^{k-1}$, is expressed as the sum of those who manage to get on train *k* and those who do not achieve it, $ee_{ij}^k$, (Equation (16)). These kinds of variables are so important when measuring quality of service because they represent passengers that are no served considering two consecutive trains after their arrival to certain stations, so it is very important to maintain their values close to zero. In the same way than before, waiting variables should be aggregated by destinations and origins, as we did above with arrivals, obtaining necessary variables like $EAD_i^k$ (Equation (17)) and $EAO_j^k$ (Equation (18)) to balance waiting people between two consecutive trains.

$$E_{ij}^k = es_{ij}^k + ee_{ij}^k, \quad \forall i,j \in H, \ i < j, \ k \in N. \tag{0}$$

$$EAD_i^k = \sum_{\{j \in H \,:\, j > i\}} E_{ij}^k, \quad \forall i \in H \setminus \{S\}, \ k \in N. \tag{0}$$

$$EAO_j^k = \sum_{\{i \in H \,:\, i < j\}} E_{ij}^k, \quad \forall j \in H \setminus \{1\}, \ k \in N. \tag{0}$$

Hence, just at the moment in which train *k* leaves the *i*-th station, the number of passengers who stay on the platform satisfy Equations (19) and (20). Both constraints type are used to balance people by origin and destination, as illustrated in Figure 4.

$$EAO_j^{k+1} = EAO_j^k + NAO_j^k - SAO_j^k, \quad \forall i \in H \setminus \{1\}, \ k \in N \setminus \{K\}. \tag{0}$$

$$EAD_i^{k+1} = EAD_i^k + NAD_i^k - SAD_i^k, \quad \forall i \in H \setminus \{S\}, \ k \in N \setminus \{K\}. \tag{0}$$

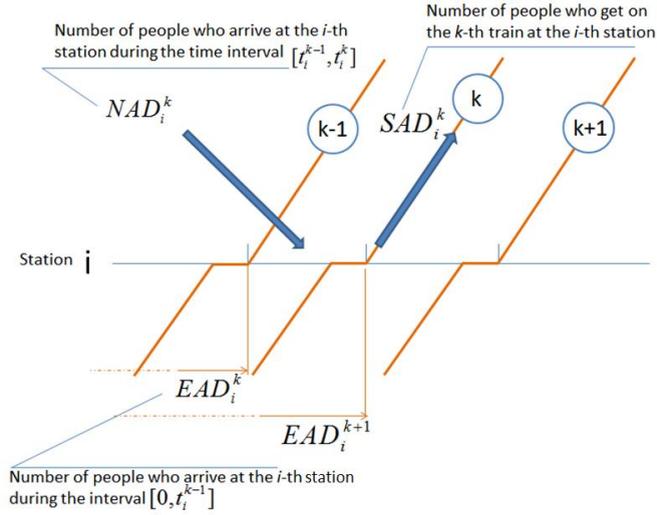

Figure 4. Waiting balance by destination

In a similar way to waiting variables, the number of passengers who get on the *k*-th train is composed of those who have arrived in the interval $\left[t^{k-1},t^{k}\right]$ and get on the *k*-th train and those who having missed the former train (($k$-1)-th train) manage to get on this *k*-th train (for each pair *i, j* of origin-destination stations), (see Equation (21)). In this case, it also becomes necessary to aggregate people by destination, $SAD_i^k$ and origin $SAO_j^k$, as it is shown in Equations (22) and (23), respectively.

$$S_{ij}^k = ns_{ij}^k + es_{ij}^k, \quad \forall i,j \in H, \ i<j, \ k \in N. \tag{0}$$

$$SAD_i^k = \sum_{\{j \in H \ : \ j >i\}} S_{ij}^{\ k}, \quad \forall i \in H \setminus \{S\}, \ k \in N. \tag{0}$$

$$SAO_j^k = \sum_{\{i \in H \ : \ i < j\}} S_{ij}^{\ k}, \quad \forall i \in H \setminus \{1\}, \ k \in N. \tag{0}$$

Now, it is possible to balance the capacity of each train at each one of the stations, using the variables defined above and the available capacity variable, $FS_i^k$. This balance is shown in Equations (24) and (25), where the last one refers to the first station.

$$FS_i^k = FS_{i-1}^k + SAO_i^k - SAD_i^k, \quad \forall i \in H \setminus \{1\}, \ k \in N. \tag{0}$$

$$FS_1^k = CAP_k - SAD_1^k, \quad k \in N. \tag{0}$$

Stop times were introduced in minimum headway constraints (Equation (7)). At this point, it is possible to determine appropriate stop intervals considering the minimum between the time needed to occupy the train available capacity (Equations (26), (27)) or the passenger flow time according to people getting on each train, Equation (28).

$$t_{stop}(i,k) \leq \text{Mint}_{stop} + g \ FS_{i-1}^k, \quad \forall k \in N, \ i \in S \setminus \{1\}. \tag{0}$$

$$t_{stop}(1,k) \leq \text{Mint}_{stop} + g \ CAP_k, \quad \forall k \in N. \tag{0}$$

$$t_{stop}(i,k) \leq \text{Mint}_{stop} + g \ SAD_i^k, \quad \forall k \in N, \ i \in S \setminus \{1\}. \tag{0}$$

$$N_{ij}^k, S_{ij}^k, E_{ij}^k, ns_{ij}^k, ne_{ij}^k, es_{ij}^k, ee_{ij}^k \geq 0, \quad \forall i,j \in H, \ i \neq j, \ k \in N.$$

$$EAD_i^k, EAO_i^k, FS_i^k, t_i^k, NAD_i^k, NAO_i^k, SAD_i^k, SAO_i^k \geq 0, \quad \forall i,j \in H, \ i \neq j, \ k \in N.$$

$$E_{ij}^1 = 0, \quad \forall i,j \in H, \ i \neq j, \quad CAP_k \geq 0, \ \forall k \in N.$$

$$\delta^k, \ \beta_i^k (0,1), \quad \forall i \in H, \ k \in N.$$

In case of a periodic schedule was preferred, it is possible to define the time between departure of two consecutive trains as the difference between departure times at the origin station, in a similar way as described in Zhou and Zhong, 2005.

$$y^k = t_1^{k+1} - t_1^k, \quad k \in N \setminus \{K\}. \tag{0}$$

We will consider, by now, a single interdeparture interval for the entire planning horizon, i.e., for the interval [0,*T*]. In a more general case, it may be more convenient consider different interdeparture times according to different planning intervals. However, we will discuss the simplest case, and we will show that in the general case of dynamic demand and capacity decisions, periodic timetables can be inappropriate, at least for some specific intervals corresponding to demand rush hours.

Zhou and Zhong, 2005 consider a minimization of the variance of interdeparture times instead of minimizing average waiting times. They justify this approach by means of the so-called Random Incidence Theorem, proposed by Larson and Odoni, 1981. Following the aforementioned theorem, the expected waiting times are proportional to the variance of interdeparture times.

$$E(w) = \frac{Var(y)}{2E(y)} + \frac{E(y)}{2}. \tag{0}$$

Where *w* refers to average waiting time and *y* represents the interdeparture time distribution. In the case of a timetable with even interdeparture times, i.e., $Var(y)=0$, then the expected waiting time is exactly half of interdeparture times. Equation (30) was obtained considering an invariable demand characterization along a specific planning interval and supposing that the vehicle capacity is large enough to accommodate all the arrivals between two consecutive trains. However, this issue is infrequent when considering rapid systems such as subways. These over-waiting times have a very important influence on the passenger's perception about quality of service and must be minimized jointly with average waiting times. To the best of our knowledge, over-waiting times have not been considered explicitly before in order to analyze quality of timetables.

With the aim of obtaining accurate timetables, we will consider a main objective-function that minimizes the total average waiting time (*AWT*) during the whole planning horizon (Equation (31)). Since it is not possible to formulate an explicit relationship between average waiting time and interdeparture times, as we commented above, it is necessary to express directly the average waiting time considering arrivals between each pair of trains and departure times at each station.

$$Min \frac{1}{\sum_{i=1}^{S-1}\sum_{j=i+1}^{S} F_{ij}^{[0,T]}} \left[ \sum_{i=1}^{S-1}\sum_{k=1}^{K}(t_i^k - t_i^{k-1})EAD_i^k + \frac{1}{2}\sum_{i=1}^{S-1}\sum_{k=1}^{K}(t_i^k - t_i^{k-1})NAD_i^k \right]. \tag{0}$$

## 3.3. Measures of timetable quality

As we mention before, the model is also used as a framework to compute measures of timetables quality and to offer to the service provider a trade-off between quality of service and operation cost measures. The importance of a multi-objective approach to analyze this kind of problems is evident. A discussion on this topic can be found in Zack, 2011. For each line and a fixed global capacity (capacity offered for the whole planning horizon in the line), the optimal *AWT* is obtained by varying parametrically the number of trains in which the global capacity is distributed, i.e., for each fixed capacity a Pareto frontier will be obtained by means of the inclusion of the corresponding constrains into the model, following the ε-constrains method proposed by Marglin, 1967. To solve this set of problems, binary variables $\delta^k$ are fixed to its corresponding values, according to the number of launched trains. Also, stopping variables are considered with value 1, i.e., each train stop at every station. Now, the model becomes a *NLP* problem, losing its integer condition.

At the next step, it would be desirable to select the most convenient capacity. With this goal, other objectives, related to the profitability of the offered service, are also considered. Specifically, horizontal, vertical and total load factors per train ($HLF_k$, $VLF_k$, $LF_k$) and served demand per train ($SD_k$) are proposed.

Load factors are, at the same time, measures of quality of service in terms of comfort (related to train occupancy) and a measure of resource usage (referred to the capacity). In this sense, the profitability of a specific service has a certain correlation with the load factor. The total load factor (commonly called load factor) aims to measure the adequacy of supply (measured in *capacity.km*) and demand (measured in *passenger.km*).

$$LF_k = \frac{\sum_{i=1}^{S-1}(Cap_k - FS_i^k)\,lng(i,i+1)}{\sum_{i=1}^{S-1}Cap_k\,lng(i,i+1)} \frac{Passenger.Km}{Capacity.Km}. \tag{0}$$

However, the study of the load factor in a service is more complex and each case must be analyzed in order to know its impact on the profitability margin. A certain value of load factor can be achieved in many different ways with diverse effects on service management. This is followed by the fact that the load factor has its origin in the combination of different phenomena, which should be identified. An important distinction comes from whether the load factor was obtained by taking many passengers during a short part of the route or few travelers during a long part of it. This distinction can be quantified by the horizontal and vertical load factors.

So, vertical load factor (*VLF*) gives an idea of the maximum train occupancy, and it would be similar, but not identical, to the above exposed concept of load factor. *VLF* is obtained by dividing the maximum passenger number ($Pm_k$) by the capacity for a certain train $k$. Note that $Pm_k$ refers to the number of passengers at the segment of the route where the higher number of passengers is achieved. Therefore, for train $k$:

$$VLF_k = \frac{Pm_k}{Cap_k} = \max_{i \in H\setminus\{S\}}\left\{\frac{Cap_k - FS_i^k}{Cap_k}\right\}. \tag{0}$$

The horizontal load factor (*HLF*) is a relation between the passenger's average trip length and the length of the line. *HLF* is obtained by dividing the average passenger number ($Pa_k$) by the maximum passenger number ($Pm_k$), above described. In some way, horizontal load factor is a measure of occupancy homogeneity along the line.

$$HLF_k = \frac{Pa_k}{Pm_k} = \left( \frac{\sum_{i=1}^{S-1}\left(Cap_k - FS_i^k\right)lng(i,i+1)}{\sum_{i=1}^{S-1}lng(i,i+1)} \right) \Big/ \max_{i \in H\setminus\{S\}}\left\{Cap_k - FS_i^k\right\}. \tag{0}$$

Moreover, it can be observed that:

$$LF_k = VLF_k \cdot HLF_k = \frac{Pm_k}{Cap_k}\frac{Pa_k}{Pm_k}$$

$$= \max_{i \in H\setminus\{S\}}\left\{\frac{Cap_k - FS_i^k}{Cap_k}\right\} \cdot \left( \left( \frac{\sum_{i=1}^{S-1}\left(Cap_k - FS_i^k\right)lng(i,i+1)}{\sum_{i=1}^{S-1}lng(i,i+1)} \right) \Big/ \max_{i \in H\setminus\{S\}}\left\{Cap_k - FS_i^k\right\} \right) \tag{0}$$

$$= \left( \frac{\sum_{i=1}^{S-1}\left(Cap_k - FS_i^k\right)lng(i,i+1)}{\sum_{i=1}^{S-1}lng(i,i+1)} \right) \Big/ Cap_k = \frac{\sum_{i=1}^{S-1}\left(Cap_k - FS_i^k\right)lng(i,i+1)}{Cap_k \sum_{i=1}^{S-1}lng(i,i+1)} = \frac{Pa_k}{Cap_k}.$$

Served demand ($SD_k$) is defined as the average percentage of people getting on train *k* in comparison to people who arrive at each station. Served demand is a measure related to the revenue, the more people getting on, the greater the service profitability and attended demand. It is important to mention that non-served demand can lead to a loss of passengers in a medium term.

$$SD_k = \frac{\sum_{\{i,j \in H: j>i\}} S_{ij}^k}{\sum_{\{i,j \in H: j>i\}} N_{ij}^k}. \tag{0}$$

Both measures, *LF* and *SD*, are relevant from the user and the service provider point of view. The confrontation of these two measures and the resulting *AWT* is done in a second phase of the decision making process in order to provide a mechanism for selecting an appropriate global capacity level and the resulting timetabling as a function of the operating cost (Lindner and Zimmermann, 2005).

The model has been implemented using *GAMS* and solved by a branch and bound procedure based on the branch and cut module *Cbc* and the cut-generation library *Cgl* included in the *COIN-OR* distribution. To solve the continuous *NLPs* the interior-point solver *Ipopt* by Watcher and Biegler, 2005, has been used. For a detailed discussion on *COIN-OR* open source library see *www.coin-or.org*.

## 4. An application of the model

To illustrate the model introduced, an application on a piece of the C5 line belonging to the Madrid suburban railway network has been tested. With this aim, we first have used surveys data (offered by RENFE) to adjust the cumulative demand functions among six stations of the line C5: Móstoles-Soto, Móstoles, Las Retamas, Alarcón, San José and Cuatro Vientos.

As we pointed out before, this adjust procedure is carried out from a minimization of the sum of quadratic errors obtained as the difference between the sigmoid approximation and the data obtained for a standard day. From now on, for the experiments shown, we assume that the carriages have a fixed capacity, equal to 40 passengers and we consider a global demand of 11581 passengers. The global capacity is increased from 2000 to 10000 passengers (2000 each time). Because of the discrete nature of the capacity, first the total capacity is divided by an integer number of wagons and then, the obtained total carriage number is distributed in one, five, ten, twenty five and fifty trains, respectively. For examples of the latter size (using the solving procedure mentioned above in an Intel Core i5 CPU 64 bits machine) typical computation time to obtain optimal solutions is about three minutes.

To perform the experiment, twenty five timetables in form of time-space diagrams (Lindner, 2000) have been obtained. On the horizontal axis of these diagrams the time is represented and on the vertical axis is shown the distance from the first station. Figure 5 shows a couple of examples corresponding to a capacity of 4000 and 8000 passengers for twenty five and fifty trains. In order to analyze the full range of results, some extreme and unrealistic values are also included (e.g., one train and 250 wagons).

Each picture shows simultaneously the number of trains, depicted as thick vertical lines, optimum departure/arrival times from/at each station as well as a part of the cumulative demand functions between each pair of stations. It can be observed how the departure times are close to the demand peaks.

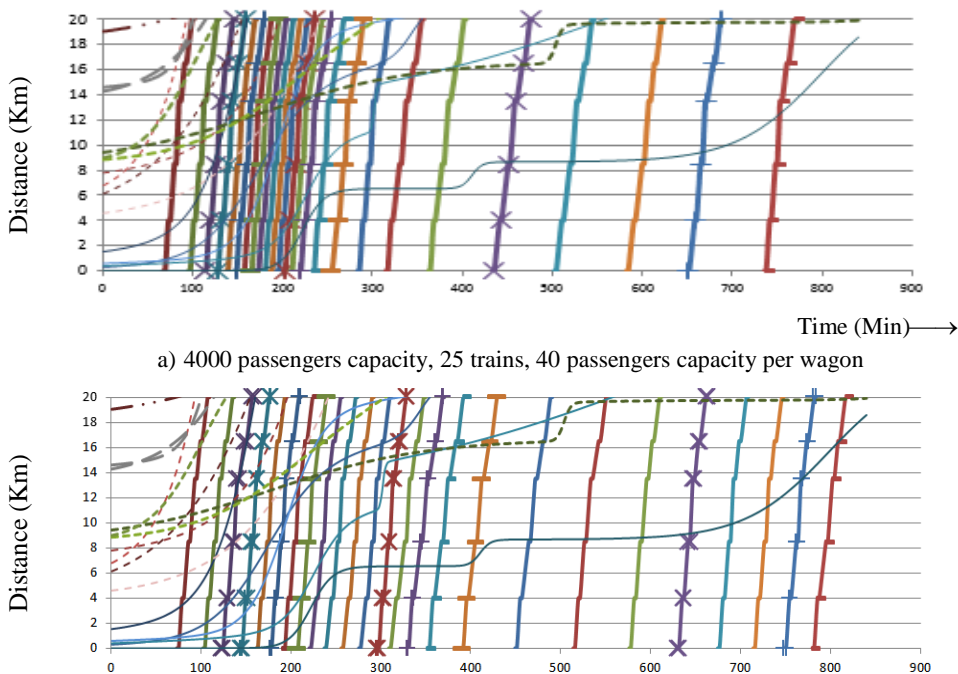

a) 4000 passengers capacity, 25 trains, 40 passengers capacity per wagon

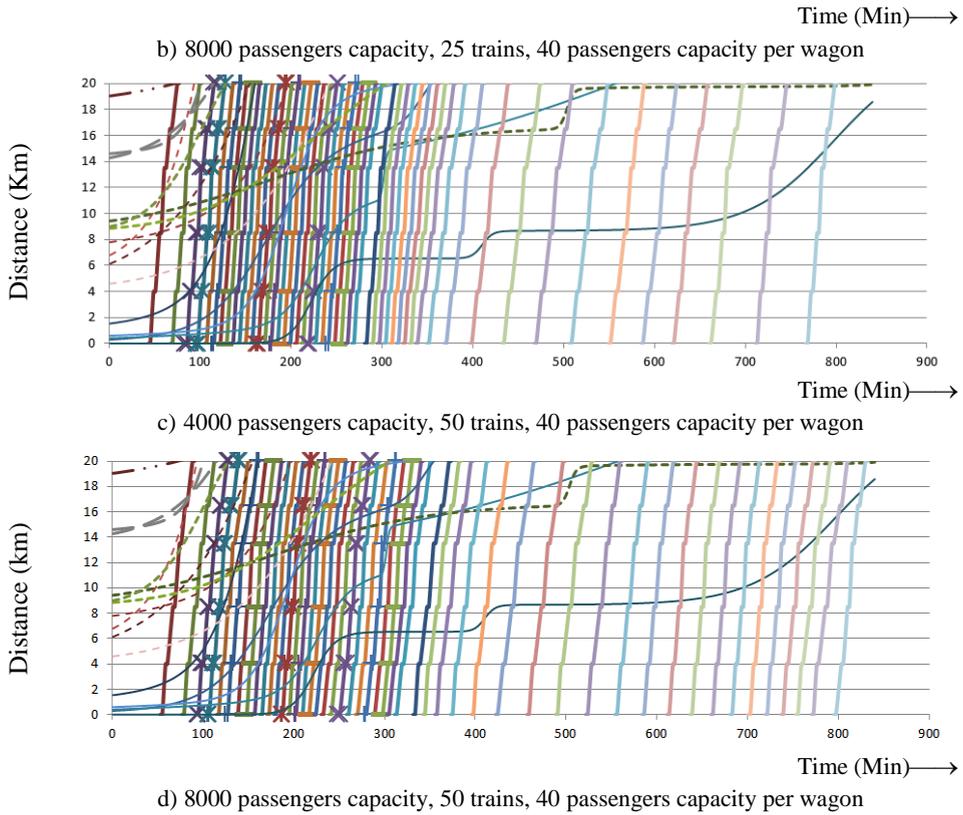

b) 8000 passengers capacity, 25 trains, 40 passengers capacity per wagon

c) 4000 passengers capacity, 50 trains, 40 passengers capacity per wagon

d) 8000 passengers capacity, 50 trains, 40 passengers capacity per wagon

Figure 5. Time-space diagrams for 25 and 50 trains, with 4000 and 8000 passengers capacity

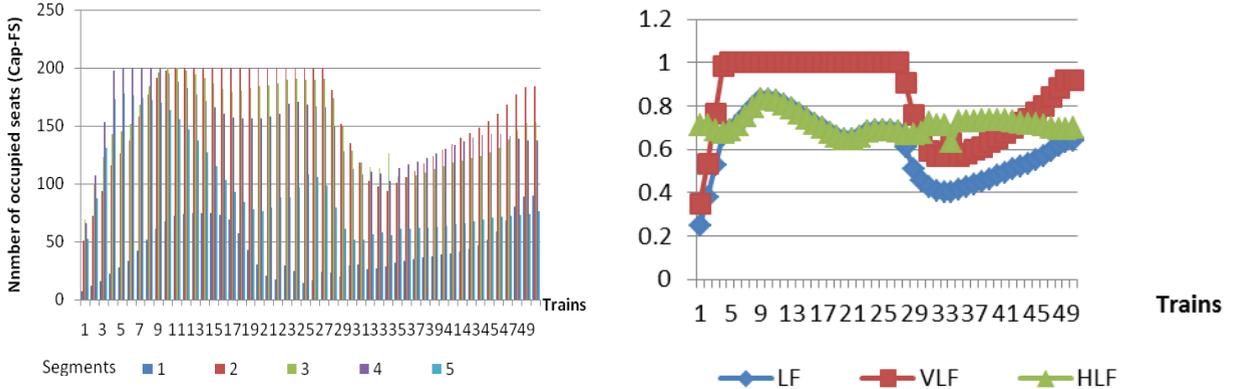

a) Number of passengers at train $i$ between stations $e_j$ and $e_{j+1}$

b) Total, Vertical and Horizontal Load Factors

Figure 6. Occupation and different Load Factors for fifty trains

Concerning to the use of train capacity, Figure 6 depicts resulting measures for the case of a total number of fifty trains launched, each one composed of five wagons. Figure 6a represents the train occupancy in each segment $l \in L$ of the line, which is directly obtained from the model solution and computed as train capacity minus the available capacity. Figure 6b shows the Total, Vertical and Horizontal Load Factors, obtained as explained in Section 3.3. The *x*-axis of both figures represents each one of the trains considered for this case. Both figures illustrate the use of the capacity of each train, however, Figure 6a does not consider the length of the segments whereas Figure 6b does it.

Some relations between both figures should be underlined, e.g., when a train reaches its full capacity during any track of its route (see Figure 6a), then, in Figure 6b, it can be observed that the *VLF* equals one and so, the *HLF* and *LF* coincide. This situation occurs between the 5$^{th}$ and 29$^{th}$ trains. A lower *VLF* (e.g., 33$^{rd}$ train) suggests that the train is not complete at any track of the route, as it can be checked on Figure 6a. This fact together with a higher value of the *HLF* represent that the available capacity is homogeneously high during the whole route. This aspect could be taken into consideration for a future redesign of the route.

The relationship between the total number of trains and the *AWT* is illustrated on Figure 7. As explained above, for each fixed capacity a Pareto frontier is obtained. In particular, Figure 7 represents the set of pairs (train number, *AWT*) for each one of the twenty five optimal timetables. The points corresponding to those timetables depicted on Figure 7 are marked by circles.

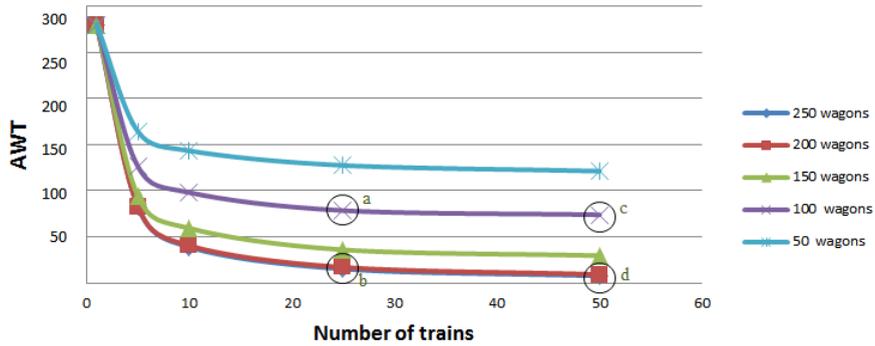

Figure 7. Pareto frontiers between Train number and *AWT* for a fixed global capacity

On the graphics, note that when the number of trains increases, then the *AWT* decreases. What is more, there exist situations in which fewer carriages, distributed among a larger number of trains, give way to a shorter *AWT*. In addition to this, differences on *AWT* between cases with 8000 and 10000 passengers capacity (200 and 250 wagons, respectively) are practically null and, since solutions with higher number of wagons are more expensive from a rolling stock point of view, in this case, it seems advisable not to consider a global capacity of 10000 passengers. All this suggests the need of a more detailed analysis of the different ways of distributing the global capacity. With this goal, each one of the obtained timetables are characterized by the Average Served Demand (*ASD*) and the Average Load Factor (*ALF*) computed from the corresponding values of $SD_k$ and $LF_k$, respectively,

$$ASD = \frac{1}{K}\sum_{k=1}^{K} SD_k, \quad ALF = \frac{1}{K}\sum_{k=1}^{K} LF_k.$$

Table 2 shows all the obtained measures. The first column indicates the total number of launched wagons during the time horizon *T*. The number of trains in which this number of wagons is distributed is indicated in the second column. The resulting above described Average Served Demand (*ASD*), Average Load Factor (*ALF*), Average Vertical Load Factor (*AVLF*) and Average Horizontal Load Factor (*AHLF*) are shown in columns third, fourth, fifth and sixth, respectively.

| Number of Wagons | Number of Trains | ASD | ALF | AVLF | AHLF |
|---|---|---|---|---|---|
| 250 | 1 | 20.67 | 11.51 | 0.24 | 1.36 |
|     | 5  | 76.06 | 48.74 | 5.41  | 6.27  |
|     | 10 | 88.97 | 55.59 | 12.77 | 12.08 |
|     | 25 | 95.64 | 59.28 | 34.49 | 29.83 |
|     | 50 | 96.88 | 59.97 | 69.89 | 59.55 |
|     | 1  | 20.67 | 14.39 | 0.29  | 1.36  |
|     | 5  | 76.04 | 60.91 | 6.76  | 6.27  |
| 200 | 10 | 87.52 | 68.49 | 15.74 | 12.07 |
|     | 25 | 92.87 | 72.11 | 41.43 | 30.21 |
|     | 50 | 93.89 | 72.80 | 83.26 | 60.72 |
|     | 1  | 20.67 | 19.18 | 0.39  | 1.36  |
|     | 5  | 68.56 | 72.59 | 7.99  | 6.32  |
| 150 | 10 | 72.92 | 76.40 | 16.67 | 12.73 |
|     | 25 | 77.07 | 78.25 | 41.67 | 32.60 |
|     | 50 | 77.43 | 78.43 | 83.33 | 65.35 |
|     | 1  | 20.67 | 28.78 | 0.59  | 1.36  |
|     | 5  | 55.01 | 80.97 | 8.33  | 6.75  |
| 100 | 10 | 57.62 | 82.52 | 16.67 | 13.75 |
|     | 25 | 59.13 | 83.31 | 41.67 | 34.71 |
|     | 50 | 59.58 | 83.51 | 83.33 | 68.76 |
|     | 1  | 20.67 | 57.55 | 1.18  | 1.36  |
|     | 5  | 39.04 | 84.42 | 8.33  | 7.03  |
| 50  | 10 | 40.69 | 86.82 | 16.67 | 14.47 |
|     | 25 | 40.82 | 87.31 | 41.67 | 36.38 |
|     | 50 | 41.06 | 89.00 | 83.33 | 74.16 |

Table 2. Characteristic measures of obtained timetables

As shown above, both measures *ALF* and *ASD* have the same behavior, i.e., both of them increase as the number of trains increases. However, a characteristic of these percentages is that they are opposed with respect to the offered capacity. Hence, for a constant number of trains, the higher the capacity, the smaller the *ALF* and larger the *ASD*.

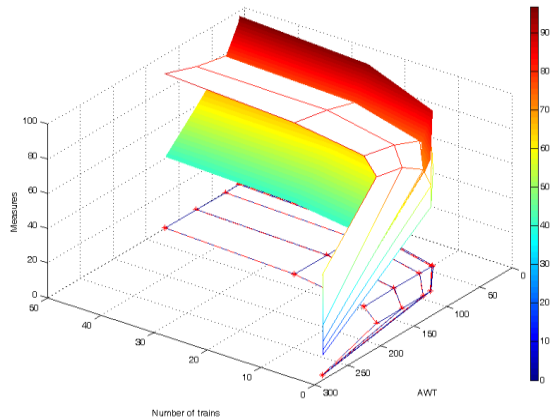

Figure 8. *ASD* and *ALF* Surfaces, level curves, *AWT* and number of trains behavior

Reaching a balance between both measures would help to find appropriate solutions. With that aim, both measures will be represented on the *z*-axis, generating the two surfaces on Figure 8. The curves representing the relationship between the total number of trains and the *AWT* are also plotted on the *xy*-plane of this figure.

Figure 9 shows the projections of both surfaces and their corresponding level curves onto *xy*-plane. Points on the intersection curve are solutions with the same *ALF* and *ASD*. Both regions *ALF > ASD* and *ALF < ASD* are depicted. The first one corresponds to better solutions from an economical point of view (a better use of the train capacity) and the second one refers to a better quality of service from user's point of view. In this experiment, taking into account that even in the best case (10000 passengers capacity) the capacity is not enough to serve all the demand along the whole day (11581 passengers), the best *AWT* are obviously near to 200-250 wagons curves and solutions with good values of *ALF* and *ASD* should be searched for values with a number of trains higher than 50. The intersection or equilibrium curve is obtained for values of the *ALF* and *ASD* around 80%; moreover, this curve is close to the Pareto frontier obtained for a total number of 150 wagons. Ideally, starting from this equilibrium between the *ALF* and *ASD*, the service provider can determine the convenient number of trains in which the total number of wagons should be distributed considering the cost per train-km (see Lindner and Zimmermann, 2005) and the penalty of the increment of the *AWT* (related to the cost per *passengers.km*).

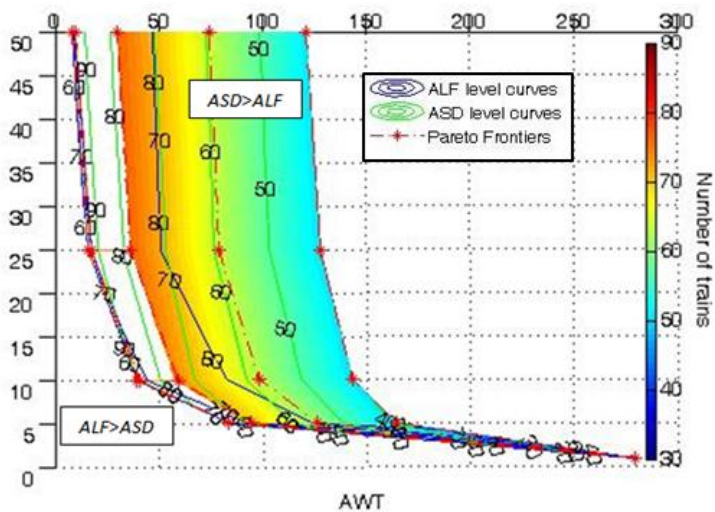

Figure 9. Projection of *ASD* and *ALF* surfaces on *xy*-plane

Another interesting possibility would be to consider, instead of fixed global capacity as in the examples above, fixed capacity per train, namely the same number of wagons, and vary the number of launched trains at each curve (see Figure 10). The obtained Pareto frontiers have the same behavior as the ones presented in Figure 7. Indeed, some values are the same, but the configuration of the curves differs in the sense that here, points of the same curve correspond to different number of trains with the same number of wagons instead of the same global number of wagons distributed into different number of trains.

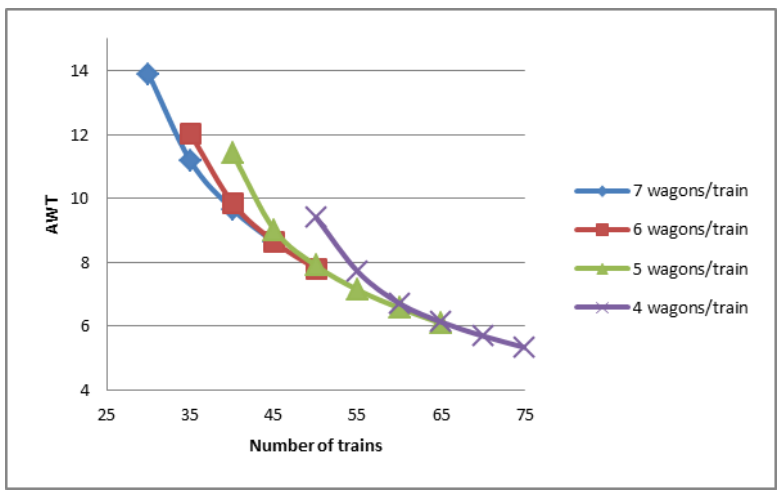

Figure 10. Pareto frontiers between Train number and *AWT* for a fixed capacity per train

From a service provider point of view, the curves in Figure 10 allow for drawing conclusions such as that not always more capacity trains provide better results like e.g., the same *AWT* is obtained when 50 trains of either five or six wagons per train are launched.

## 5. Conclusions

This paper introduces the determination of non-cyclic railway timetables relaxing the unrealistic assumptions of regular demand behavior and fixed capacity. A non-linear integer programming model is presented to decide how many trains should be scheduled in order to attend a given variable demand and compute the departure and stop times of each scheduled train for all their tour stops. The model is also used as a basis to obtain measures related to the quality of timetables attending quality of service and line profitability, being a rather different and improved version of the preliminary methodology described in Canca et al., 2011.

Number of services, train capacity, stops, speeds and arrival/departure times are optimized simultaneously, starting with a description of the demand between each pair of stations of a given line. A continuous demand representation in terms of a continuous origin-destination matrix for the full day operation has been used instead of a set of discrete *O-D* demand matrices corresponding to different time intervals or the peak hour *O-D* matrix. Each matrix element represents the evolution of the daily travel demand between each pair of stations. In order to incorporate the demand description into the model, the cumulative demand curve is used. In fact, cumulative demand is approximated by a linear combination of a certain number of sigmoidal functions which improves the one considered in Canca et al., 2011, because it ensures a good approximation outside the planning horizon. This allows us to produce a more accurate demand description. The parameters of each term are determined using a least square procedure.

Optimal timetables are computed adjusting arrival and departure times to the above described cumulative demand functions by means of passenger arrival behavior with the objective of minimizing the average waiting time. Some constraints to incorporate aspects like periodicity and to determine the optimal number of trains have been also discussed.

In order to validate the model, several tests have been done and it can be checked that the obtained measures correspond with the expected ones. So, the frequency of trains is higher close to the demand peaks and the relation between the load factor and the usual occupancy can be easily observed in the graphics. Moreover, the flexibility of the model allows for the selection of different objectives, as well as a post-study of the obtained results by varying parametrically relevant inputs such as the number of trains or their capacity.

Timetables quality can be analyzed following different points of view. Frequently, optimization models offer the best timetable under one or two criteria, commonly minimizing average waiting time (considering a uniform demand and representing the user's point of view) or minimizing differences of departure times with respect to an ideal timetable designed by the infrastructure manager (service provider point of view). But there are other aspects that must be considered like capacity or vehicles use. The proposed model allows for the evaluation of performance measures of timetables considering arrival and departure times as input data. This feature enables us to compute several ratios to characterize the quality of timetables considering additional issues like train capacity and occupancy.

With the aim of reducing the search field for appropriate solutions as well as discarding non-efficient solutions, we compute a set of Pareto frontiers for the bi-objective problem which takes into account the number of trains and the passenger average waiting time (*AWT*). These trade-offs are obtained following the ε-constrains method and fixing parametrically the global capacity of the line. Then, the Average Served Demand (*ASD*), which represents the average number of people who can get on each train at any station, and the Average Load Factor (*ALF*), related to the profitability of each train service, are introduced in order to define a new trade-off that depict the equilibrium between occupancy and attended demand. This equilibrium curve, in conjunction with the *AWT*, allows the service provider to take a decision on the adequate number of trains and the most suitable train capacity guaranteeing a certain level of service quality.

This paper emphasizes the need of a deep analysis of the timetables obtained with different optimization models considering alternative measures concerning quality of service and line profitability and presents a methodology to develop such kind of studies. A final analysis of the network performance for each line would require the inclusion of different aspects like maintenance, rolling stock and crew management.

## 6. Acknowledgments

This research work has been partially supported by the Excellence program of the Andalusian Government, under grant P09-TEP-5022 and AZ under grant P11-FQM-7276. Moreover, we are grateful to RENFE (Madrid suburban railway network) for providing data and information relevant to the realistic adjustment of the model. A.Z. also acknowledges financial support from Ministerio de Educación y Ciencia of Spain under grants MTM2009–14668–C02–02 and MTM2011–28952–C02, cofinanced by the European Community fund FEDER and from Universidad Politécnica de Madrid, Spain ("TACA" consolidated research group).